\documentclass[a4paper,10pt]{article}
\usepackage{amscd}

\usepackage{amsthm}

\newtheoremstyle{pedro}{}{}{\itshape}{}{\sc}{~--}{ }{\thmname{#1}\thmnumber{ #2}\thmnote{ (#3)}}

\newtheoremstyle{pedroremark}{}{}{\rm}{}{\sc}{~--}{ }{\thmname{#1}\thmnumber{ #2}\thmnote{ (#3)}}

\theoremstyle{pedro}
\newtheorem{lem}{Lemma}[subsection]

\newtheorem{thm}[lem]{Theorem}
\newtheorem*{thm2}{Theorem}

\newtheorem{prop}[lem]{Proposition}

\newtheorem{coro}[lem]{Corollary}

\theoremstyle{pedroremark}

\newtheorem{rmk}[lem]{Remark}

\newtheorem{notations}[lem]{Notations}

\theoremstyle{definition}

\newtheorem{ex}[lem]{Example}

\usepackage{amsmath, amssymb, amsfonts, amscd}
\usepackage{amsthm}

\newcommand{\Z}{\mathbf{Z}}
\newcommand{\C}{\mathbf{C}}
\newcommand{\Spin}{\mathbf{Spin}}
\renewcommand{\O}{\mathbf{O}}
\newcommand{\SO}{\mathbf{SO}}

\newcommand{\GL}{\mathbf{GL}}
\newcommand{\SL}{\mathbf{SL}}
\newcommand{\Sp}{\mathbf{Sp}}
\newcommand{\Gm}{\mathbf{G}_m}

\newcommand{\aff}{\mathbf{A}}
\newcommand{\pr}{\mathbf{P}}
\newcommand{\spec}{\mathrm{Spec}}
\newcommand{\F}{\mathbf{F}}
\newcommand{\bor}{_}
\newcommand{\et}{\mathrm{et}}
\newcommand{\Mu}{\boldsymbol{\mu}}

\usepackage{titlesec}

\titleformat{\section}{\bf\center}{\S \arabic{section}.}{0.5em}{}
\titleformat{\subsection}{\bf}{\arabic{section}.\arabic{subsection}.}{.5em}{}
\newcommand{\lc}{{\em loc. \!\!\!\! cit.}}
\title{Geometric methods for cohomological invariants}
\author{Pierre Guillot}

\begin{document}

\maketitle

\begin{abstract}
We explain how to exploit Rost's theory of Chow groups with
coefficients to carry some computations of cohomological invariants. In
particular, we use the idea of the ``stratification method''
introduced by Vezzosi.

We recover a number of known results, with very different proofs. We
obtain some new information on spin groups.

\bigskip

{\bf Status:} This paper has been accepted for publication in 
{\em Documenta Mathematica}.

\end{abstract}

%\tableofcontents

\section{Introduction}

{\em In what follows, $k$ is a base field, and all other fields
considered in this paper will be assumed to contain $k$. We fix a
prime number $p$ which is different from $char(k)$, and for any field
$K/k$, we write $H^i(K)$ for $H^i(Gal(K_s/K), \Z/p(i))$. Here $K_s$
is a separable closure of $K$, and $\Z/p(i)$ is the $i$-th Tate twist
of $\Z/p$, that is $\Z/p\otimes\Mu_p^{\otimes i}$. (See the end of this introduction for more information on our choice of coefficients.) }

\subsection{Cohomological invariants}

Given a functor $A(-)$ from fields over $k$ to the category of pointed
sets, the (mod $p$) {\em cohomological invariants} of $A$ are all
transformations of functors
$$a: A(-) \to H^*(-).$$ Typical examples for $A$ include
$$\begin{array}{rcl} 
A(K) & = & \textrm{isomorphism classes of nondegenerate quadratic forms over $K$},\\ 
A(K) & = & \textrm{isomorphism classes of octonions algebras over $K$},\\
\textrm{etc...} & &
\end{array}$$

For more examples and a very good introduction to the subject, see
\cite{serre}. 

We are mostly interested in the situation where $A(K)$ is the set of isomorphism classes over K of some type of "algebraic object", especially
when all such objects become isomorphic over algebraically closed fields. In this case, if we write $X_K$ for the base point in $A(K)$, then we obtain a group scheme $G$ over $k$ by setting $G(K) = Aut_K(X_K)$, assuming $Aut_K$ is appropriately defined. Moreover, by associating to each object $Y_K\in A(K)$ the variety of isomorphisms (in a suitable sense) from $Y_K$ to the base point $X_K$, we obtain a $1:1$ correspondence between $A(K)$ and
$H^1(K,G)$. Recall that this is the (pointed) set of all $G$-torsors
over $K$, i.e. varieties acted on by $G$ which become isomorphic to
$G$ with its translation action on itself when the scalars are
extended to the algebraic closure $\bar K$. 

For example, in the case of quadratic forms, resp. of octonion
algebras, we have $A(K)=H^1(K,\O_n)$, resp. $A(K)=H^1(K,G_2)$.

In more geometric terms, a $G$-torsor over $K$ is a principal
$G$-bundle over $\spec (K)$. In this way we see that cohomological
invariants $$H^1(-,G) \to H^*(-)$$ are analogous to characteristic
classes in topology. We shall call them the mod $p$ cohomological
invariants of $G$. They form a graded $H^*(k)$-algebra which we denote
by $Inv(G)$ (the prime $p$ being implicit).

\subsection{Versal torsors and classifying spaces}\label{sub:class}
A simple way of constructing a $G$-torsor is to start with a
$G$-principal bundle $E\to X$ and take the fibre over a point $\spec
(K)\to X$. As it turns out, there always exist bundles $\pi:E\to X$ with
the particularly nice feature that {\em any} $G$-torsor is obtained,
up to isomorphism, as a fibre of $\pi$. In this case, the generic
torsor $T_\kappa$ over $\kappa=k(X)$ is called {\em versal}. We shall
also say that $\pi$ is a versal bundle.

A strong result of Rost (presented in \cite{serre}) asserts that an
invariant $a\in Inv(G)$ is entirely determined by its value
$a(T_\kappa)$ on a versal torsor. It follows that $Inv(G)\subset
H^*(\kappa)$, and more precisely one can show that $$Inv(G)\subset
A^0(X,H^*).$$ The right hand side refers to the cohomology classes in
$H^*(\kappa)$ which are {\em unramified} at all divisors of $X$ (see
\lc). (The notation comes from Rost's theory of Chow groups with coefficients, see below and section \ref{sec:rost}.)

There are two well-known constructions of versal bundles. One
can embed $G$ in a ``special'' group $S$ (eg $S=\GL_n$, $S=\SL_n$ or
$\Sp_n$) and take $E=S$, $X=S/G$. Alternatively, one can pick a
representation $V$ of $G$ such that the action is free on a nonempty
open subset $U\subset V$, and take $E=U$, $X=U/G$.

In either situation, there are favorable cases when we actually have
$$Inv(G)=A^0(X,H^*).$$ For versal bundles of the first type, this
happens when $S$ is simply connected (see Merkurjev \cite{merkur}). For
versal bundles of the second type, this happens when the complement of
$U$ has codimension $\ge 2$ (\cite{serre}, letter by Totaro).

In this paper we shall restrict our attention exclusively to the
second type of versal torsors, for reasons which we shall explain in a
second. We shall view $U/G$ as an approximation to the {\em classifying
space} of $G$ (see \cite{totaro}). It will not cause any confusion, hopefully, to call this variety $BG$ (see \S\ref{subsec:equiv}).

\subsection{Rost's Chow groups with coefficients}
The group $A^0(X,H^*)$, for any variety $X$, is the first in a
sequence of groups $A^n(X,H^*)$, Rost's {\em Chow groups with
  coefficients}. They even form a (bigraded) ring when $X$ is
smooth. The term ``Chow groups'' is used since the usual mod $p$ Chow
groups of $X$ may be recovered as $CH^n X \otimes_\Z \Z/p=A^n(X,H^0)$.

Rost's groups have extremely good properties: homotopy invariance,
long exact sequence associated to an open subset, and so on. There is
even a spectral sequence for fibrations.

Our aim in this paper, very briefly, is to use these geometric
properties in order to get at $A^0(BG,H^*)=Inv(G)$.

To achieve this, we shall be led to introduce the equivariant Chow
groups $A^n_G(X,H^*)$, defined when $X$ is acted on by the algebraic
group $G$. The definition of these is in perfect analogy with the
case of equivariant, ordinary Chow groups as in \cite{edidin},
\cite{totaro}, where ``Borel constructions'' are used. When $X=\spec
(k)$, we have $A^n_G(\spec (k),H^*)=A^n(BG,H^*)$ where $BG$ is a
classifying space ``of the second type'' as above.

The key observation for us will be that $A^0_G(V,H^*)=Inv(G)$ when $V$
is a representation of $G$ (thus a variety which is equivariantly
homotopic to $\spec (k)$, so to speak). This allows us to cut out $V$
into smaller pieces, do some geometry, and eventually implement the
``stratification method'', which was first introduced by Vezzosi in
\cite{vezzosi} in the context of ordinary Chow groups.

\subsection{Results}

A large number of the results that we shall obtain in this paper are
already known, although we provide a completely different approach for
these. Occasionally we refine the results and sometimes we even get
something new.

In this introduction we may quote the following theorem (\ref{thm:exactseq} in the text):

\begin{thm2} Let $p=2$ and $n\ge 2$. There are exact sequences:
$$\begin{CD}
0 @>>> Inv(\O_n) @>>> Inv(\Z/2 \times \O_{n-1}) @>r>> Inv(\O_{n-2})
\\0 @>>> Inv(\SO_{2n}) @>>> Inv(\Z/2 \times \SO_{2n-1}) @>r>> Inv(\SO_{2n-2})
\\0 @>>> Inv(\Spin_{2n}) @>>> Inv(\Z/2 \times \Spin_{2n-1}) @>r>> Inv(\Spin_{2n-2})
\end{CD}$$

Moreover, the image of the map $r$ contains the image of the restriction map $Inv(\O_{n})\to Inv(\O_{n-2})$, resp. $Inv(\SO_{2n})\to Inv(\SO_{2n-2})$, resp. $Inv(\Spin_{2n})\to Inv(\Spin_{2n-2})$.
\end{thm2}

In each of the three cases the existence of the map $r$ was not known,
as far as I am aware, while the other half of the exact sequence is
described by Garibaldi in \cite{skip}. Of course the computation of
$Inv(\O_n)$ and $Inv(\SO_n)$ has been completed (see \cite{serre}
again), and indeed we shall derive it from the theorem. This affords a new construction of an invariant first defined by Serre (and written originally $b_1$). On the other hand the invariants of $\Spin_n$ are not known in general: they have been computed for $n\le 12$ only (\cite{skip}), so our theorem might be of some help.

% Another result (\ref{thm:pgl4} below) which is perhaps worth
% mentioning concerns $\PGL_4$. Before stating it, let us recall that the invariants of $\PGL_n$ when $n$ is not a prime are considerably harder to deal with than the mod $p$
% invariants of $\PGL_p$ (see \cite{skip}). Nonzero invariants of $\PGL_4$ have been shown to exist for $p=2$ in degree $2$ and $4$ so far (see \cite{tignol}), where we say that an invariant has degree $d$ when it takes its values in $H^d(-)$. It is rumoured that Rost knows how to prove that there are no invariants of degree $>4$ (other than $H^*(k)$-linear combinations of invariants of lower degree). This result however does not appear in the literature (or anywhere, as far as I am aware), so our weaker result may still be interesting:
% 
% \begin{thm2}
% Let $p=2$ and let $k$ be algebraically closed. When $d>11$, we have
% $$Inv^d(\PGL_4)=0.$$\end{thm2}
% 
% In this statement, and throughout this paper, we write $Inv^d(G)$ for
% the invariants of $G$ of degree $d$. 

% In any case, our emphasis in this paper is with methods rather than specific results. For example, the invariants of $\PGL_n$ are usually worked out using the theory of simple, central algebras (see the introduction to section \ref{sec:pgln}). By contrast, we rely on the representation theory of $\PGL_n$, which seems new.

In any case, our emphasis in this paper is with methods rather than specific results, and the point that we are trying to make is that the stratification method is a powerful one. It provides us with a place to start when trying to tackle the computation of $Inv(G)$, whatever $G$ may be. It is much more mechanical than any other approach that I am aware of.

Still, the reader will find in what follows a number of examples of computations of $Inv(G)$ for various $G$'s: products of copies of $\Mu_p$ in \ref{prop:elem}, the group $\Spin_7$ in \ref{ex:spin7}, the wreath product $\Gm\wr \Z/2$ in \ref{ex:kleenex}, the dihedral group $D_8$ in \ref{ex:d8}, etc. 

\subsection{Organization of the paper}

We start by presenting Rost's definition of Chow groups with
coefficients in section \ref{sec:rost}. We also indicate how to
construct the equivariant Chow groups, and we mention a number of
basic tools such as the K\"unneth formula. 

In section \ref{sec:strat}
we introduce the stratification method, and prove the above results on
orthogonal groups. We compute $Inv(\O_n)$ and $Inv(\SO_n)$ completely.

In section \ref{sec:proj} we explain how one could use the
stratification method with projective representations rather than
ordinary ones. We recover a result (corollary \ref{coro:skip}) which
was proved and exploited fruitfully in \cite{skip}.

We proceed to introduce the Bloch \& Ogus spectral sequence in section
\ref{sec:spec}. %This will complete our review of methods. The spectral sequence is well-known, and so is its relationship to cohomological invariants, but it is convenient to mention at this point what it gives with equivariant groups and how it fits with known calculations of ordinary Chow rings of classifying spaces. 
As we have mentioned already, the stratification method has been used by many authors (including, Vezzosi, Vistoli and myself) to compute $CH^* BG$ and $H^*BG$, and to study the cycle map $CH^* BG \to H^* BG$. It turns out that the spectral sequence shows, roughly speaking, that cohomological invariants measure the failure of this cycle map to be an isomorphism (in small degrees this is strictly true). This was the motivation to try and extend the stratification method to cohomological invariants.

%Next we prove in section \ref{sec:pgln} some general results on the invariants of $\PGL_n$, and we specialize them to $n=2$ and $n=4$. This part is rather technical.

We conclude in section \ref{sec:milnor} with some easy remarks on invariants with values in other cycle modules, particularly algebraic $K$-theory.

\smallskip
\noindent{\em Notations \& Conventions.} We insist on the assumptions
that we have made at the beginning of this introduction: $p$ is a
fixed prime, $k$ is a fixed base field, $p\ne char (k)$, and $H^*(K)$
means mod $p$ Galois cohomology as defined more precisely above. Our particular choice of twisting for the coefficients has been dictated by the desire to obtain a "cycle module" in the sense described in section \ref{sec:rost}, and at the same time keep things as simple as possible. There are more general cycle modules, including Galois cohomology simply with some more general coefficients (see \cite{rost}). In most applications, we have $p=2$ anyway, in which case the twisting does not interfere.

We will say that $G$ is an algebraic group to indicate that $G$ is a
smooth, affine group scheme over $k$ [in alternative terminology: $G$
is a linear algebraic group over $\bar k$ defined over $k$]. Subgroups
will always be assumed to be closed and smooth. We shall encounter
many times the stabilizer of an element under an action of $G$: in
each case, it will be trivially true that the scheme-theoretic
stablizer is closed and smooth. [Alternatively: subgroups are always
defined over $k$. So will all stabilizers we encounter.]

We write $Inv(G)$ for the mod $p$ cohomological invariants of $G$,
even though the letter $p$ does not appear in the notation. We caution
that in \cite{skip}, these would be called the invariants of $G$ with
values in $\Mu_p$ (for $p=2$ there is no difference).

\smallskip
\noindent{\em Acknowledgements.} I would like to thank Burt Totaro and Skip Garibaldi for their many suggestions.

\section{Rost's Chow groups with coefficients}\label{sec:rost}

\subsection{Definitions}\label{subsec:chowrost}

In \cite{rost}, Rost defines a {\em cycle module} to be a functor $M:
F \mapsto M(F)$ from fields containing $k$ to graded abelian groups,
equipped with structural data. Quoting Rost, these divide into "the
even ones", namely restriction and corestriction maps, and "the odd
ones": it is required that $M(F)$ should have the structure of a
module over $K_*(F)$ (Milnor's $K$-theory of $F$), and there should
exist residue maps for discrete valuations. The fundamental example
for us is $M(F)=H^*(F)$ (mod $p$ Galois cohomology as defined in the
introduction). Another example is Milnor's $K$-theory itself,
$M(F)=K_*(F)$.

Given a cycle module $M$, Rost defines for any variety $X$ over $k$
the Chow groups with coefficients in $M$, written $A_i(X,M)$ ($i\ge
0)$. These are bigraded; for the cycle modules above, we have the
summands $A_i(X,H^j)$ and $A_i(X,K_j)$. Moreover, the classical and
mod $p$ Chow groups can be recovered as $CH_iX=A_i(X,K_0)$ and $CH_iX
\otimes_\Z \Z/p=A_i(X,H^0)$.

When $X$ is of pure dimension $n$, we shall reindex the Chow groups by
putting $A^i(X,H^j)=A_{n-i}(X,H^j)$. Our interest in the theory of
Chow groups with coefficients stems from the concrete description of
$A^0(X,H^j)$ when $X$ is irreducible: it is the subgroup of
$H^j(k(X))$ comprised of those cohomology classes in degree $j$ which
are unramified at all divisors of $X$. Hence, when $X$ is a
classifying space for an algebraic group $G$ (as in \ref{sub:class}),
we have $A^0(X,H^j)=Inv^j(G)$.

Rost's Chow groups have all the properties of ordinary Chow groups,
and more: in fact they were designed to be more flexible than ordinary
Chow groups, particularly in fibred situations. For the time being, we
shall be content to list the following list of properties; we will
introduce more as we go along. Here $X$ will denote an equidimensional
variety.

\begin{enumerate}
\item When $X$ is smooth, $A^*(X,H^*)$ is a graded-commutative ring.
\item A map $f: X\to Y$ induces $f^*: A^*(Y,H^*)\to A^*(X,H^*)$
whenever $Y$ is smooth or $f$ is flat.
\item A proper map $f: X\to Y$ induces $f_*: A_*(X,H^*)\to A_*(Y,H^*)$.
\item There is a projection formula: $f_*(xf^*(y))=f_*(x)y$ (here
$f:X\to Y$ is proper, and $X$ and $Y$ are smooth).
\item Let $i:C\to X$ be the inclusion of a closed subvariety, and let
$j:U\to X$ denote the inclusion of the open complement of $C$. Then
there is a long exact sequence:
\[ \mspace{-100mu} \begin{CD} \cdots @>>> A_m(C,H^n) @>{i_*}>> A_m(X,H^n) @>{j^*}>>
A_m(U,H^n) @>{r}>> A_{m-1}(C,H^{n-1}) @>>> \cdots
\end{CD}\]  
Moreover the residue map $r$, or connecting homomorphism, satisfies $$r(xj^*(y))=r(x)i^*(y).$$
\item If $\pi:E\to X$ is the projection map of a vector bundle or an
affine bundle, then $\pi^*$ is an isomorphism.
\end{enumerate}

There are particular cases when we can say more about the exact
sequence in property (5). When $n=0$ for instance, we recover the
usual localisation sequence:
$$\begin{CD}
CH_mC \otimes_\Z \Z/p @>>> CH_mX\otimes_\Z \Z/p @>>> CH_mU\otimes_\Z \Z/p @>>> 0.
  \end{CD}$$

Also, when $X$ is equidimensional and $C$ has codimension $\ge 1$ (for
example when $X$ is irreducible and $U$ is nonempty), then the
following portion of the sequence is exact:
$$\begin{CD}
0 @>>> A^0(X,H^n) @>{j^*}>> A^0(U,H^n).
  \end{CD}$$

\begin{ex}\label{ex:gm}
Let us compute $A^0(\Gm,H^*)$, where $\Gm=\GL_1$ is the punctured
affine line. The inclusion of the origin in $\aff^1$ yields the exact
sequence:
$$\mspace{-50mu}\begin{CD}
0 @>>> A^0(\aff^1,H^*)  @>>> A^0(\Gm,H^*) @>>> A^0(\spec (k), H^{*-1}) @>>> A^1(\aff^1, H^{*-1}).
  \end{CD}$$
  
Now, $\aff^1$ can be seen as a (trivial) vector bundle over $\spec
(k)$, so $A^i(\aff^1,H^*)=A^i(\spec (k), H^*)$. From the definitions,
we see directly that $A^i(\spec (k), H^*)=0$ when $i>0$, and
$A^0(\spec (k), H^*)=H^*(k)$.

So we have the exact sequence of $H^*(k)$-modules:
$$\begin{CD}
0 @>>> H^*(k) @>>> A^0(\Gm,H^*) @>>> H^{*-1}(k) @>>> 0.
  \end{CD}$$ 
It follows that $A^0(\Gm,H^*)$ is a free $H^*(k)$-module on two
generators, one in degree $0$, the other in degree $1$. In fact, if we
write $k(\Gm)=k(t)$, it is not difficult to see that the element $(t)
\in H^1(k(t))=k(t)^\times / (k(t)^\times)^p$ can be taken as the
generator in degree $1$. For this one only has to unwind the
definition of $r$ as a residue map: see \cite{serre}, chap. II and
III. Essentially the reason is that the divisor of the function $t$ on
$\aff^1$ is the origin.
  
  In any case, we shall use the letter $t$ for this generator. When
  $k$ is algebraically closed, so that the ring $H^*(k)$ can be seen
  as the field $\F_p$ concentrated in degree $0$, then we can write
  $A^0(\Gm,H^*)=\F_p[t]/(t^2)$, an exterior algebra.
\end{ex}
  
\begin{rmk}\label{rmk:norm}
Let $X$ be smooth, and suppose that $X$ has a $k$-rational point
$\spec (k) \to X$. The induced map on Chow groups gives a splitting
for the map $H^*(k) \to A^0(X,H^*)$ coming from the projection $X \to
\spec (k)$. It follows that we can write $$A^0(X,H^*)=H^*(k) \oplus
A^0(X,H^*)_{norm}.$$ The elements in the second summand are said to be
{\em normalised}. Note that the splitting may {\em a priori} depend on
the choice of a $k$-rational "base point" for $X$.

There is an analogous notion for cohomological invariants. Namely, an
invariant $n\in Inv(G)$ is called normalised when it is $0$ on the
trivial torsor; an invariant $c\in Inv(G)$ is called constant when
there is an $\alpha\in H^*(k)$ such that for any torsor $T$ over a
field $K$, the value $c(T)$ is the image of $\alpha$ under the natural
map $H^*(k)\to H^*(K)$. Any invariant $a\in Inv(G)$ may be written
$a=n + c$ with $n$ normalised and $c$ constant, and we may write
$$Inv(G)=H^*(k) \oplus Inv(G)_{norm}$$ where the second summand
consists of normalised elements.

Whenever $X$ is a classifying space of the types considered in the
introduction, there is a canonical choice of a point $\spec (k) \to X$
such that the pullback of the versal torsor over $X$ is the trivial
torsor. It follows that $A^0(X,H^*)_{norm}=Inv(G)_{norm}$.

Let us give an application. The group $\Mu_p$ acts freely on $\Gm$, and
the quotient is again $\Gm$. Since $\Gm$ is special, it follows from
the introduction that the torsor $\Gm\to\Gm$ is versal, so
$Inv(\Mu_p)\subset A^0(\Gm,H^*)$. It is not immediate that this
inclusion is an equality, as $\Gm$ is not a classifying space: $\Gm$
is not simply-connected, and the origin has only codimension $1$ in
$\aff^1$, so we are in neither of the two favorable cases.

To show that it is indeed an equality, we remark that
$A^0(\Gm,H^*)_{norm}$ is a free $H^*(k)$-module on one generator of
degree $1$, from the previous example. As a result, if we can find a
non-zero, normalised invariant of $\Mu_p$ in degree $1$, then the map
of $H^*(k)$-modules $Inv(\Mu_p)_{norm} \to A^0(\Gm,H^*)_{norm}$ will
certainly be surjective. Such an invariant is easy to find. Indeed,
take the identity $H^1(K,\Mu_p) \to H^1(K)$. (Should you try to compute
the mod $p$ invariants of $\Mu_\ell$, where $\ell$ is a prime
different from $p$, you would find no normalized invariants at all. It
is an easy exercise to prove this using the techniques to be
introduced in this paper).
\end{rmk}

\subsection{Spectral sequences and K\"unneth formula}\label{subsec:kunneth}

Apart from Milnor K-theory and Galois cohomology, there is one extra
type of cycle module which we shall consider in this paper. In the
presence of a map $f: X \to Y$, and having picked a primary cycle
module $M$, there is for each $n\ge 0$ a new cycle module written
$A_n[f;M]$, see \cite{rost},\S 7. These are functors from fields
$\kappa$ "over $Y$", that is fields with a map $\spec (\kappa) \to Y$,
to graded abelian groups -- which turns out to be enough to define
$A_*(Y,A_n[f;M])$. Quite simply, if we write $X_\kappa=\spec (\kappa)
\times_Y X$, the definition is:
$$A_n[f;M](\kappa)=A_n(X_\kappa,M).$$

When $f$ and $M$ are understood, we shall write $\mathcal{A}_n$ for
$A_n[f;M]$. Also, when $f$ is flat, so that the fibres all have the
same dimension $d$, we may change the grading to follow the
codimension, thus defining $\mathcal{A}^n=\mathcal{A}_{d-n}$.

\begin{thm}\label{thm:specseq} 
Let $f: X\to Y$ be a flat map, and suppose that $Y$ is
equidimensional. Let $\mathcal{A}_n=A_n[f;M]$ for a cycle module
$M$. Then there exists a convergent spectral sequence:
$$A^r(Y,\mathcal{A}^s)\Rightarrow A^{r+s}(X,M).$$
\end{thm}
  
See \cite{rost},\S 8, for a proof.
  
\begin{coro}\label{coro:kunneth}
Let $X$ be a scheme (over $k$), and write $X_\kappa$ for the scheme
obtained by extending the scalars to the field $\kappa$. Assume that
$A^0(X,H^*)$ is a free $H^*(k)$-module of finite rank, and that
$$A^0(X_\kappa,H^*)=A^0(X,H^*) \otimes_{H^*(k)} H^*(\kappa).$$
Then there is a K\"unneth isomorphism:
$$A^0(X\times Y,H^*)=A^0(X,H^*)\otimes_{H^*(k)} A^0(Y,H^*)$$ for any
equidimensional $Y$.
\end{coro}

The hypotheses of this corollary should be compared with \cite{serre}, 16.5.

\begin{proof} 
Let $f: X \times Y \to Y$ be the projection. By the theorem, we know
that $A^0(X\times Y, H^*)=A^0(Y, \mathcal{A}^0)$.

However, for any field $\spec (\kappa)\to Y$, the fibre of $f$ above
$\spec (k)$ is $X_\kappa$ as in the statement of the corollary. Our
hypothesis is then that $\mathcal{A}^0(\kappa)=A^0(X_\kappa, H^*)$ is
a direct sum of copies of $H^*(\kappa)$, and more precisely that
$\mathcal{A}^0$ splits up as the direct sum of several copies of the
cycle module $H^*$. The result follows.
\end{proof}

As an illustration, we have the following proposition.

\begin{prop}\label{prop:elem} There are invariants $t^{(i)}\in Inv((\Mu_p)^n)$ for
  $1\le i \le n$, each of degree $1$, such that $Inv((\Mu_p)^n)$ is a
  free $H^*(k)$-module on the products $t^{(i_1)}t^{(i_2)}\cdots
  t^{(i_k)}$, for all sequences $1\le i_1 < i_2 < \cdots <i_k \le n$.

When $k$ is algebraically closed, $Inv((\Z/p)^n)$ is an exterior
algebra $\Lambda(t^{(1)},\ldots,t^{(n)})$ over $\F_p$.
\end{prop}

\begin{proof} 
The case $n=1$ has been dealt with in remark \ref{rmk:norm}. The
general case follows from corollary \ref{coro:kunneth}.\end{proof}

%The next corollary generalises property (6) in \S
%\ref{subsec:chowrost}.
  
%\begin{coro}[of theorem \ref{thm:specseq}]\label{coro:fibresaffine} Suppose $f:E\to X$ is a map such that for each point $\spec (\kappa) \to X$, the fibre $\spec (\kappa) \times_X E$ is isomorphic to an affine space (as a $\kappa$-variety). Then 
%$$A^*(E,H^*)=A^*(X,H^*).$$ 
%\end{coro}

%\begin{proof} Using the spectral sequence, we reduce the problem immediately to the case $X=\spec (\kappa)$. Then $E=\aff^n_\kappa$. As did Rost in order to prove property (6), we note that an induction reduces further the problem to $n=1$, in which case the desired isomorphism exists in virtue of the axioms of cycle modules.
%\end{proof}

\subsection{Equivariant Chow groups}\label{subsec:equiv}

The good properties of Chow groups with coefficients will allow us to
define the equivariant Chow groups of a variety acted on by an
algebraic group. This will be in strict analogy with \cite{edidin} and
\cite{totaro}, to which we will refer for details.

So let $G$ be a linear algebraic group over $k$. There exists a
representation $V$ of $G$ such that the action is free outside of a
closed subset $S$; moreover, we can make our choices so that the
codimension of $S$ is arbitrarily large (\lc). Put $U=V-S$.

Let $X$ be any equidimensional scheme over $k$ with an action of
$G$. Define
\[X \bor G=\frac{U\times X}{G}.\] 
We shall restrict our attention to those pairs $(X,U)$ for which
$X\bor G$ is a scheme rather than just an algebraic space. Lemma 9 and
Proposition 23 in \cite{edidin} provide simple conditions on $X$ under
which an appropriate $U$ may be found. This will be amply sufficient
for the examples that we need to study in this paper.\footnote{I am
  grateful to the referee for pointing out this technical difficulty. I
  also agree with him or her that it would be desirable to extend
  Rost's theory of Chow groups with coefficients to algebraic
  spaces. This is certainly not the place to do so.}

The notation $X\bor G$ hides the dependence on $U$ because, as it
turns out, the Chow groups $A^i(X \bor G,H^*)$ do not depend on the
choice of $V$ or $S$, as long as the codimension of $S$ is large
enough (depending on $i$). One may prove this using the "double
fibration argument" as in \cite{totaro}. We write $A^i_G(X,H^*)$ for
this group (or the limit taken over all good pairs $(V,S)$, if you
want).

A map $f:X\to Y$ induces a map $f_G: X \bor G \to Y \bor G$ (the
notation $X \bor G$ will mean that a choice of $V$ and $S$ has been
made). It follows that the equivariant Chow groups with coefficients
$A^*_G(-,H^*)$ are functorial, and indeed they have all the properties
listed in \S\ref{subsec:chowrost}. The proof of this uses that if $f$
is flat, proper, an open immersion, a vector bundle projection, etc,
then $f_G$ is respectively flat, proper, an open immsersion, or a
vector bundle projection (see \cite{edidin}). It may be useful to
spell out that a $G$-invariant open subset $U$ in an equidimensional
$G$-variety $X$, whose complement $C$ has codimension $c$, gives rise
to a long exact sequence:

$$
\cdots \to A^m_G(C,H^n) \to A^{m+c}_G(X,H^n) \to A^{m+c}_G(U,H^n) \to A^{m+1}_G(C,H^{n-1}) \to \cdots
.$$

We shall refer to it as the equivariant long exact sequence associated to $U$.

\begin{ex}
 We shall write $BG$ for $\spec (k) \bor G$. It follows from the
 results above and from the introduction that $$A^0_G(\spec
 (k),H^*)=A^0(BG,H^*)=Inv(G).$$
\end{ex}

\begin{ex}\label{ex:freeactions}
 Suppose that the action of $G$ on $X$ is free, and that the quotient
 $X/G$ exists as a scheme. Then there is a natural map $X\bor G \to
 X/G$. Moreover $X\bor G$ is an open subset in $(V\times X)/G$, which in
 turn is (the total space of) a vector bundle over $X/G$. The
 complement of $X\bor G$ in $(V\times X)/G$ can be taken to have an
 arbitrarily large codimension (namely, it is that of $S$). In a given
 degree $i$, we may thus choose $V$ and $S$ appropriately and obtain:
$$A^i_G(X,H^*)=A^i(X\bor G,H^*)=A^i((V\times X)/G,H^*)=A^i(X/G,H^*).$$
(The second equality coming from the long exact sequence as in
\S\ref{subsec:chowrost}, property (5).)

Thus when the action is free, the equivariant Chow groups are just the Chow groups of the quotient.
\end{ex}

\begin{notations}
As we have already done in this section, we shall write $X\bor G$ to signify that we have chosen $V$ and $S$ with the codimension of $S$ large enough for the computation at hand. Thus we can and will write indifferently $A^*_G(X,H^*)$ or $A^*(X\bor G,H^*)$.

Moreover, we shall write $EG$ for $U=V-S$, with the same convention. The quotient $EG/G$ is written $BG$, and we write $A^*(BG,H^*)$ much more often than $A^*_G(\spec (k),H^*)$.
\end{notations}

\begin{rmk}\label{rmk:complex}
Whenever we have a theory at hand which has the properties listed in \S\ref{subsec:chowrost}, we can define the equivariant analog exactly as above. Apart from Chow groups with coefficients, we shall have to consider at one point the \'etale cohomology of schemes, namely the groups $H^i_\et(X,\Z/p(i))$, which we shall write simply $H^i_\et(X)$, the coefficients being understood. In fact we shall only encounter the expression $H^i_\et(BG)$.

When $k=\C$, we have $H^*_\et(BG)=H^*(BG,\F_p)$, where on the right hand side we use topological cohomology and a model for the classifying space $BG$ in the classical, topological sense of the word. To see this, note first that $H^*_\et((V-S)/G)=H^*((V-S)/G,\F_p)$ since we are using finite coefficients. Moreover, we can arrange $V$ and $S$ so that $V-S$ is a Stiefel variety (see \cite{totaro}), and this provides sufficiently many $V$'s and $S$'s (that is, the codimension of $S$ can be arbitrarily large, with $V-S$ a Stiefel variety). As a result, we can restrict attention to these particular pairs $(V,S)$ when forming the limit $H^*_\et(BG)=\lim H^*((V-S)/G,\F_p)$ without changing the result, which is then clearly $H^*(E/G,\F_p)$ where $E$ is the infinite Stiefel variety. This space $E$ is contractible, so $E/G$ is a topological model for $BG$.
\end{rmk}

\section{The stratification method}\label{sec:strat}

\subsection{The idea of the method}
\label{subsec:strat}
Let $G$ be a linear algebraic group. The stratification method is a
procedure to compute $A^0(BG,H^*)=Inv(G)$. We shall not try to present
it as a mechanical algorithm, but rather as a heuristic recipe. This
being said, one of the virtues of the method is that it is closer to a
systematic treatment of the question than any other approach that the
author is aware of.

The stratification method rests on two very elementary facts.

\begin{enumerate}

\item If $K$ is a closed subgroup of $G$, then the equivariant Chow
groups of the variety $G/K$ have an easy description: indeed there is an isomorphism
of varieties:
$$(G/K)\bor G = \frac{EG \times G/K}{G} = \frac{EG}{K}=BK.$$
Thus $A^*_G(G/K,H^*)=A^*(BK,H^*).$

\item Suppose that $V$ is a representation of $G$. Then it is a
$G$-equivariant vector bundle over a point, and therefore $V\bor G \to
\spec (k)\bor G=BG$ is a vector bundle, too. As a result,
$A^*_G(V,H^*)=A^*(BG,H^*)$.
\end{enumerate}

Things are put together in the following way. One starts with a
well-chosen representation $V$, and then cuts $V$ into smaller
pieces. These smaller pieces would ideally be orbits or families of
orbits parametrized in a simple manner. Using (1), and hoping that the
stabilizer groups $K$ showing up are easy enough to understand, one
computes the Chow groups of the small pieces. Applying then repeatedly
the equivariant long exact sequence from \S\ref{subsec:equiv}, one
gets at $A^0_G(V,H^*)$. By (2), this is $A^0(BG,H^*)$.

Any particular application of the method will require the use of {\em
ad hoc} geometric arguments. Before giving an example however, we
propose to add a couple of lemmas to our toolkit.

\begin{lem}\label{lem:biactions}
Suppose that $G$ and $K$ are algebraic groups acting on the left on
$X$, and suppose that the actions commute. Assume further that the
action of $K$ alone gives a $K$-principal bundle $X \to X/K=Y$. Then
$$A^*_{G}(Y,H^*)=A^*_{G\times K}(X,H^*).$$
\end{lem}

\begin{proof}
This is simply an elaboration on example \ref{ex:freeactions}.
\end{proof}

We shall have very often the occasion of using this lemma in the
following guise:

\begin{lem}\label{lem:family} Suppose that $G$ contains a subgroup $G'$ which is an extension 
$$1\to N\to G' \to K\to 1.$$ Let $G\times K$ act on $G/N$ by
$(\sigma,\tau)\cdot [g]=[\sigma g \tau^{-1}]$. Assume moreover that
$K$ acts on a variety $X$, and that there is a $K$-principal bundle
$G/N \times X \to Y.$

Then $$A^*_G(Y,H^*)=A^*_{G'}(X,H^*).$$ If the product $N\times K$ is
direct, we have $X\bor N \times K=BN \times X\bor K$.
\end{lem}

Here the variety $Y$ is an example of what we call "a family of orbits
(of $G$) parametrized in a simple manner". Note that when $K$ is the
trivial group and $X=\spec (k)$, we recover property (1) above.

\begin{proof} By lemma \ref{lem:biactions}, we have
$$A^*_G(Y,H^*)=A^*_{G\times K}(G/N \times X, H^*).$$

Now let $E_1$ be an open subset in a representation of $G$, on which
the action is free, and suppose that the codimension of the complement
of $E_1$ is large enough. Similarly, pick $E_2$ for $K$; then $E_1
\times E_2$ can be chosen for the group $G\times K$. We regard $E_1$
as a trivial $K$-space and $E_2$ as a trivial $G$-space. Finally, $G'$
has natural maps to both $G$ and $K$, and we combine these to see $E_1
\times E_2$ as a $G'$-space, with non-trivial action on each factor.

We write:
$$\begin{array}{rcl}
G/N \times X \bor G\times K & = & (E_1\times E_2 \times \frac{G}{N} \times X)/G\times K
\\                          & = & (\frac{E_1 \times G/N}{G} \times E_2 \times X)/K
\end{array}.$$

By property (1) above, we may identify $\frac{E_1 \times G/N}{G}$ with
$E_1/N=BN$. Moreover we have arranged things so that, under this
identification, the action of $K$ on $\frac{E_1 \times G/N}{G}$
translates into its natural action on $E_1/N$ (that is, the action
induced from that which $G'$ possesses as a subgroup of $G$
normalising $N$). We proceed:
$$\begin{array}{rcl}
G/N \times X \bor G\times K & = & (E_1/N \times E_2 \times X)/ K
\\                          & = & (E_1 \times E_2 \times X)/G'
\end{array}.$$
When the product is direct, the action of $K$ on $E_1/N$ is trivial. This concludes the proof.\end{proof}

As an example, we may apply this to actions of $G$ on disjoint unions:

\begin{lem}\label{lem:disjoint}
Suppose that $G$ acts on $Y$ and that $Y$ is the disjoint union of varieties $Y_1 \coprod \ldots \coprod Y_n$. Let $X=Y_1$, and suppose that $X$ contains a point in each orbit, that is, suppose that the natural map
$$f: G\times X \to Y $$ is surjective. Let $K$ be the subgroup of $G$ leaving $X$ invariant. Then there is an identification: $$A^0_G(Y,H^*)=A^0_K(X,H^*).$$
\end{lem}

\begin{proof} This is immediate from the previous lemma, with $N$ the trivial group.
\end{proof}

\subsection{The example of $\O_n$, $\SO_n$, and $\Spin_n$}
For this example we take $p=2$ (so that $char(k)\ne 2$).

Let $q$ be a non-degenerate quadratic form on a $k$-vector space $V$ of dimension $n$. Assume moreover that $q$ has maximal Witt index. Then its automorphism group, which we will denote by $\O_n$, is split. Similarly one has the groups $\SO_n$ and $\Spin_n$, which are also split.

\begin{thm}\label{thm:exactseq} Let $p=2$ and $n\ge 2$. There are exact sequences:
$$\begin{CD}
0 @>>> Inv(\O_n) @>>> Inv(\Z/2 \times \O_{n-1}) @>r>> Inv(\O_{n-2})
\\0 @>>> Inv(\SO_{2n}) @>>> Inv(\Z/2 \times \SO_{2n-1}) @>r>> Inv(\SO_{2n-2})
\\0 @>>> Inv(\Spin_{2n}) @>>> Inv(\Z/2 \times \Spin_{2n-1}) @>r>> Inv(\Spin_{2n-2})
\end{CD}$$

Moreover, the image of the map $r$ contains the image of the restriction map $Inv(\O_{n})\to Inv(\O_{n-2})$, resp. $Inv(\SO_{2n})\to Inv(\SO_{2n-2})$, resp. $Inv(\Spin_{2n})\to Inv(\Spin_{2n-2})$.

\end{thm}
\begin{proof}
\noindent{\em Step 1.}
Let $G_n$ denote either of $\O_n$, $\SO_n$, or $\Spin_n$. Then $G_n$ has the canonical representation $V$. Let $V'$ denote $V-\{0\}$. We write $U$ for the open subset in $V'$ on which $q$ is non-zero, and we write $C$ for its complement.

The codimension of $\{ 0 \}$ in $V$ is $n$. Consider the long exact sequence:
$$\mspace{-50mu}
0 \to A^0_{G_n}(V,H^*) \to A^0_{G_n}(V',H^*) \to 0=A^{1-n}_{G_n}(\{ 0 \},H^{*-1}) \to A^1_{G_n}(V,H^{*-1}) \to A^1_{G_n}(V',H^{*-1}).
$$
It follows that $A^0_{G_n}(V,H^*) = A^0_{G_n}(V',H^*)=Inv({G_n})$. Moreover the last map above is surjective when $*=1$ (cf \S\ref{subsec:chowrost}), so that $A^1_{G_n}(V,H^0) = A^1_{G_n}(V',H^0)=CH^1B{G_n}\otimes_\Z \Z/2$.

Turning to the equivariant long exact sequence associated to the open set $U$ in $V'$, we finally get:
$$\begin{CD}
0 @>>> Inv({G_n}) @>>> A^0_{G_n}(U,H^*) @>>> A^0_{G_n}(C,H^{*-1}).\qquad (\dagger)
\end{CD}$$
When $*=1$ we have:
$$\mspace{-50mu} 0 \to Inv^1({G_n}) \to A^0_{G_n}(U,H^1) \to A^0(C,H^0) \to CH^1B{G_n} \otimes_\Z \Z/2 \to CH^1_{G_n} U\otimes_\Z \Z/2 \to 0. \quad (\dagger\dagger)$$

\smallskip
\noindent{\em Step 2.} Let $Q=q^{-1}(1)$. We shall use the following result.

\begin{lem} \label{lem:actions}
The action of $G_n$ on $Q$ is transitive. Moreover the stabilizer of a $k$-rational point is isomorphic to $G_{n-1}$, and we get an isomorphism of $k$-varieties $Q=G_n/G_{n-1}$.

The action of $G_n$ on $C$ is transitive. Moreover the stabilizer of a
$k$-rational point is isomorphic to a semi-direct product $H\rtimes
G_{n-2}$, where $H$ is an algebraic group isomorphic to affine space
as a variety, and we get an isomorphism of $k$-varieties
$C=G_n/H\rtimes G_{n-2}$. Finally, the map $B(H\rtimes G_{n-2}) \to
BG_{n-2}$ is an affine bundle.
\end{lem}

For a proof, see \cite{vistolimolina}.

From this we can at once identify the last term in the exact sequence
($\dagger$). Indeed, if $S=H\rtimes G_{n-2}$ as in the lemma, then
$A^0_{G_n}(C,H^*)=A^0(BS,H^*)$ as explained at the beginning of this
\S. However, since $BS\to BG_{n-2}$ is an affine bundle, we draw from
\S\ref{subsec:chowrost}, property (6), that
$A^0(BS,H^*)=A^0(BG_{n-2},H^*)=Inv(G_{n-2})$.

\smallskip
\noindent{\em Step 3.} We turn to the term $A^0_{G_n}(U,H^*)$. Extend
the action of $G_n$ on $Q$ to an action on $Q\times \Gm$ which is
trivial on the second factor. The group $\Z/2=\{1,\tau\}$ also acts on
$Q\times \Gm$ by $\tau (x,t)=(-x,-t)$. The two actions
commute. Finally, there is a map $Q\times \Gm \to U$ defined by
$(x,t)\mapsto tx$. It is $G_n$-equivariant, and a $\Z/2$-principal
bundle. 

We wish to apply lemma \ref{lem:family} with $N=G_{n-1}$, $K=\Z/2$,
and $G'=N\times K$. This is of course the time when we need to assume
that $n$ is even if $G_n$ is $\SO_n$ or $\Spin_n$. Then the element
$-Id$ generates the copy of $\Z/2$ that we need. We conclude:
$$A^0_{G_n}(U,H^*)=A^0(BG_{n-1} \times (\Gm)\bor{\Z/2},H^*).$$

%The action of $\Z/2$ being free, we conclude as in example
%\ref{ex:freeactions} that $A^0_{G_n}(U,H^*)=A^0_{\Z/2 \times
%G_n}(Q\times \Gm,H^*)$.
 %In this case, $G_n$ has an element acting as $-Id$ on $V$, and as a result, the group $\Z/2 \times G_n$ has the same orbits on $Q$ as $G_n$.

%Using this and the fact that the action of $G_n$ on the factor $\Gm$ is trivial, we have
%$$\begin{array}{rcl}
%(Q\times \Gm)\bor \Z/2 \times G_n & = & \frac{EG\times E\Z/2 \times Q \times \Gm}{\Z/2 \times G_n}
%\\                                 & = & \left( \frac{EG \times Q}{G_n} \times E\Z/2 \times \Gm \right) / (\Z/2)
%\\                                 & = & \frac{EG \times Q}{G_n} \times \frac{E\Z/2 \times \Gm}{\Z/2}
%\\                                 & = & (Q\bor G_n) \times (\Gm\bor \Z/2)
%\end{array}$$

We may easily compute $A^0((\Gm) \bor{\Z/2},H^*)$. Indeed, using that the
action is free, with quotient $\Gm$, this is $A^0(\Gm,H^*)$ from
example \ref{ex:freeactions}. As observed in remark \ref{rmk:norm},
this happens to equal $A^0(B\Z/2,H^*)=Inv(\Z/2)$.

Now, by the K\"unneth formula (corollary \ref{coro:kunneth}), which we
may use as $Inv(\Z/2)$ is a free $H^*(k)$-module, we draw 
$$A^0(BG_{n-1} \times (\Gm)\bor{\Z/2},H^*)=Inv(G_{n-1}) \otimes_{H^*(k)}
Inv(\Z/2)=Inv(\Z/2\times G_{n-1}).$$

\smallskip
\noindent{\em Step 4.} We need to prove the last sentence in the
theorem. We claim that there is an $x\in Inv(\Z/2 \times G_{n-1})$
such that $r(x)=1 \in Inv(G_{n-2})$. Granted this, the theorem follows
from the formula $r(xj^*(y))=r(x)i^*(y)=i^*(y)$ (see
\S\ref{subsec:chowrost}), applied here with $i^*$ denoting the
restriction $Inv(G_n)\to Inv(G_{n-2})$ and $j^*$ denoting the
restriction $Inv(G_n) \to Inv(\Z/2 \times G_{n-1})$.

To prove the claim we start by computing $CH^1_{G_n}U$ (in the rest of
this proof we shall write simply $CH$ for the mod $2$ Chow groups). As
above this is $CH^1 BG_{n-1} \times (\Gm)\bor{\Z/2}$. Arguing as in
example \ref{ex:freeactions}, we see that we may replace
$(\Gm)\bor{\Z/2}$ by $\Gm/\Z/2=\Gm$, as far as computing the Chow
groups goes. Now, $\Gm$ being an open set in affine space, it has a
trivial Chow ring and there is a K\"unneth formula, so:
$$ CH^1 BG_{n-1} \times (\Gm)\bor{\Z/2}=CH^1BG_{n-1}.$$

Therefore the last map in ($\dagger\dagger$) is the restriction map
$CH^1BG_n \to CH^1BG_{n-1}$. When $G_n=\Spin_n$, which is
simply-connected, we have $CH^1BG_n=0$. In \cite{vistolimolina}, the
reader will find a proof that $CH^1B\SO_{2n}=0$ also and that
$CH^1B\O_n \to CH^1B\O_{n-1}$ is an isomorphism. The claim follows
from the exact sequence ($\dagger\dagger$).\end{proof}

\subsection{End of the computation for $\O_n$}\label{subsec:ortho}

The main novelty in theorem \ref{thm:exactseq} is with the spin groups
(specifically, with the second half of the exact sequence presented,
as the first half was obtained by Garibaldi in \cite{skip}). The
invariants of $\O_n$ and $\SO_n$ are completely known, see
\cite{serre}. In this section and the next, however, we show how to
recover these results from theorem \ref{thm:exactseq}. In the case of
$\O_n$ this is so close to the proof given in \lc that we shall omit
most of the details. In the case of $\SO_n$ on the other hand, our
method is considerably different, so we have thought it worthwhile to
present it. We continue with $p=2$.

A major role will be played by {\em Stiefel-Whitney classes}. Recall
that $H^1(K,\O_n)$ is the set of isomorphism classes of non-degenerate
quadratic forms on $K^n$. Given such a form $q$, the classes
$w_i(q)\in H^*(K)$, for $i\ge 0$, have been defined by Milnor
\cite{milnor} (originally they were defined in mod $2$ Milnor
K-theory). They are $0$ for $i>n$ and if $w_t(q)=\sum w_i(q)t^i$, then
one has
$$w_t(q\oplus q')=w_t(q)w_t(q').$$

Each $w_i$ can be seen as a cohomological invariant in $Inv(\O_n)$. In fact one has

\begin{prop}\label{prop:on} The $H^*(k)$-module $Inv(\O_n)$ is free with a basis given by the classes $w_i$, for $0\le i \le n$. Moreover, the following multiplicative formula holds:
$$\begin{array}{rcll}
w_1w_{i-1} & = & w_i & (i~ \textrm{odd}),
\\     & = & (-1)w_{i-1} & (i~ \textrm{even}).
\end{array}$$
\end{prop}

In the formula, $(-1)$ stands for the image of $-1\in k$ in $H^1(k)=k^\times/(k^\times)^2$.

\begin{proof}[Sketch proof] From theorem \ref{thm:exactseq}, we have an injection
$$0\to Inv(\O_n) \to Inv(\Z/2 \times \O_{n-1})$$ where $\Z/2 \times \O_{n-1}$ is seen as a subgroup of $\O_n$, the nonzero element in the copy of $\Z/2$ corresponding to the $-Id$ matrix. An immediate induction then shows that $Inv(\O_n)$ injects into $Inv((\Z/2)^n)$, the invariants of the elementary abelian $2$-group comprised of the diagonal matrices in $\O_n$ with $\pm 1$ as entries. (This could be seen as a consequence of the surjection $H^1(K,(\Z/2)^n) \to H^1(K,\O_n)$, which in turns expresses the fact that quadratic forms may be diagonalized).

We have computed $Inv(\Z/2)$ in remark \ref{rmk:norm}, and we may apply corollary \ref{coro:kunneth} to obtain $Inv((\Z/2)^n)$. The symmetric group $S_n$ acts on this ring, and the image of $Inv(\O_n)$ must certainly lie in the subring fixed by this action. The latter is easily seen to be a free $H^*(k)$-module on the images of the Stiefel-Whitney classes.

The multiplicative formula is obtained by direct computation.\end{proof}

\subsection{End of the computation for $\SO_n$}

Since $\SO_n$ is a subgroup of $\O_n$, the restriction map $Inv(\O_n) \to Inv(\SO_n)$ defines Stiefel-Whitney classes for $\SO_n$.

Now, as we have defined it, the group $\SO_n$ has the following
interpretation for its torsors: $H^1(K,\SO_n)$ is the set of
isomorphism classes of non-degenerate quadratic forms on $K^n$ whose
discriminant is $1$. Recall that the discriminant is the determinant
of any matrix representing the bilinear form associated to $q$, viewed
as an element of $K^\times/ (K^\times)^2$. As it turns out, the
discriminant of $q$ is precisely $w_1(q)$ when seen as a class in the
additive group $H^1(K)$. So we have $w_1=0$ in $Inv(\SO_n)$. From the
formula $w_1w_{2i}=w_{2i + 1}$, we have in fact $w_{2i+1}=0$ in
$Inv(\SO_n)$.

\begin{prop} When $n=2m+1$ is odd, $Inv(\SO_n)$ is a free $H^*(k)$-module with basis $\{w_0, w_2, w_4, \ldots, w_{2m}\}$.

When $n=2m$ is even, $Inv(\SO_n)$ is a free $H^*(k)$-module with basis $$\{w_0, w_2, w_4, \ldots, w_{2m -2},b^{(1)}_{n-1}\},$$ where $b^{(1)}_{n-1}$ is an invariant of degree $n-1$.
\end{prop}

The invariant $b^{(1)}_{n-1}$ was introduced by Serre \cite{serre}. It was originally denoted by $b_1$, but we try to keep lowerscripts for the degree whenever possible.

\begin{proof} The odd case is easy, for we have an isomorphism
$$\Z/2 \times \SO_{2m + 1} = \O_{2m + 1}.$$ From the K\"unneth formula (corollary \ref{coro:kunneth}), it follows that $Inv(\SO_{2m+1})$ is the quotient of $Inv(\O_{2m+1})$ by the ideal generated by $w_1$. From the formula in proposition \ref{prop:on}, this is the submodule generated by the odd Stiefel-Whitney classes. The result follows.

Turning to the even case, we use the exact sequence from theorem \ref{thm:exactseq}:
$$\begin{CD}
0 @>>> Inv(\SO_{2m}) @>>> Inv(\Z/2 \times \SO_{2m-1}) @>r>> Inv(\SO_{2m-2})
\end{CD}$$

We proceed by induction, assuming the result for $n-2$.

The group $\Z/2 \times \SO_{2m-1}$ is seen as a subgroup of
$\SO_{2m}$, were again the copy of $\Z/2$ is identified with $\{\pm Id
\}$. However, there is also a canonical isomorphism of $\Z/2 \times
\SO_{2m - 1}$ with $\O_{2m-1}$, as above. The corresponding injective
map $\O_{2m-1} \to \SO_{2m}$ thus obtained induces the map
$H^1(K,\O_{2m-1}) \to H^1(K,\SO_{2m})$ which sends a quadratic form
$q$ to $\tilde q=(q\otimes \det (q)) \oplus \det (q)$. Here $\det (q)$
is the $1$-dimensional quadratic form with corresponding $1\times 1$
matrix given by the discriminant of $q$.

We have $w_{2i+1}(\tilde q)=0$ as explained above. To compute the even
Stiefel-Whitney classes, we note that $w_i(\tilde q)=w_i((q\oplus
\langle 1 \rangle)\otimes \det (q))=w_i(q\otimes \det (q))$, and we recall that if
we factorise formally
$$\sum w_i(q) t^i=\prod (1 + a_it^i),$$ then
$$\sum w_i(q\otimes \det (q)) t^i=\prod (1 + (a_i+ w_1(q))t^i).$$

It follows that $w_{2i}(\tilde q)=w_{2i}(q) + w_1(q)R_i(q)$, and the remainder $w_1(q)R_i(q)$ can thus be written as a linear combination over $H^*(k)$ of the classes $w_j(q)$ with $j<2i$. Moreover $w_{2m}(\tilde q)=0$.

Now, write $w_i$, resp $w_i'$, resp $w_i''$, for the Stiefel-Whitney
classes in $Inv(\SO_{2m})$, resp $Inv(\O_{2m-1})$, resp
$Inv(\SO_{2m-2})$, and regard $Inv(\SO_{2m})$ as a subring of
$Inv(\O_{2m-1})$. Thus we have $w_{2i+1}=0$ and $w_{2i}=w_{2i}' +
w_1'R_i$. 

Let $M$, resp $N$, denote the $H^*(k)$-submodule of
$Inv(\O_{2m-1})$ generated by the classes $w_{2i+1}'$ for $0\le i \le
m-2$, resp by the classes $w_{2i}$ for $0\le i\le m-1$ together with
$b^{(1)}_{n-1}=w_{2m+1}' + (-1)w_1'R_m$. Then $M$ and $N$ are free
modules, and $Inv(\O_{2m-1})=M\oplus N$. We claim (i) that the residue
map $r$ is $0$ on $N$ (ie $r(b^{(1)}_{n-1})=0$), and (ii) that $r$
maps $M$ injectively onto a free submodule in $Inv(\SO_{2m-2})$. This
will complete the induction step, as we will have $Inv(\SO_{2m})=\ker
(r)=N$.

Both parts of the claim are proved by the same computation. We first note from theorem \ref{thm:exactseq} that $r(w_1')=1$, and we compute: $$
\begin{array}{rcl}
r(w_{2i+1}') & = & r(w_1'w_{2i}')
\\             & = & r(w_1'w_{2i}) - r((w_1')^2R_i')
\\             & = & w_{2i}'' - r((-1)w_1'R_i)
\\             & = & w_{2i}'' - (-1)r(w_1'R_i)
\end{array}. $$ Take $i=m-1$ to obtain (i) (since $w_{2m-2}''=0$). Property (ii) is clear.

It remains to start the induction. The group $\SO_2$ is a torus, so $Inv(\SO_2)=H^*(k)$. One gets the result for $\SO_4$ from this and theorem \ref{thm:exactseq}. This is very similar to the induction step (but easier), and will be left to the reader.\end{proof}

\section{Projective variants}\label{sec:proj}
\subsection{Projective bundles}

Consider the $n$-th projective space $\pr^n$. It has an open subset isomorphic to $\aff^n$, and from the long exact sequence of \S\ref{subsec:chowrost}, (5), we draw at once $A^0(\pr^n,H^*)=H^*(k)$.

Now, if $V\to X$ is a vector bundle, we may form the associated
projective bundle $\pi:\pr (V) \to X$. The fibre of $\pi$ over $\spec
(\kappa) \to X$ is $\pr^{n-1}_\kappa$, where $n$ is the rank of
$V$. Thus we see that the induced cycle module
$\mathcal{A}^0=A^0[\pi,H^*]$ on $X$ is isomorphic to $H^*$, and from
theorem \ref{thm:specseq} we have:

\begin{lem}\label{lem:proj} Let $\pi:\pr (V) \to X$ denote the projective bundle associated to the vector bundle $V\to X$. Then
$$A^0(\pr (V),H^*)=A^0(X,H^*).$$
\end{lem}

\begin{rmk} It should be kept in mind that the above argument, which rests on the spectral sequence of theorem \ref{thm:specseq}, has been used for the sake of concision only. The result is a consequence of the following general statement:
$$A^*(\pr (V),H^*)=A^*(X,H^*)[\zeta]/(\zeta^{n})$$ with $\zeta \in A^1(\pr (V),H^0)$. This is perfectly analogous to the usual statement for ordinary Chow groups, and is no harder to prove. We shall have no use for the complete statement in the sequel, however, and therefore we omit the lengthy argument.
\end{rmk}

\begin{coro}\label{coro:skip} Let $V$ be a representation of the linear algebraic group $G$. Form the projective representation $\pr (V)$, and assume that there is a $k$-rational point in $\pr (V)$ whose orbit is open, and isomorphic to $G/S$. Then there is an injection
$$0 \to Inv(G) \to Inv(S).$$ %When this open orbit is the whole of $\pr (V)$, the injection is an isomorphism.
\end{coro}

\begin{proof} We have a vector bundle $V\bor G \to BG$ and the associated projective bundle is $\pr (V) \bor G \to BG$. From the lemma, $$Inv(G)=A^0(BG,H^*)=A^0(\pr (V) \bor G,H^*)=A^0_G(\pr (V),H^*).$$

We may view $G/S$ as an open subset in $\pr (V)$, and from the equivariant long exact sequence we have
$$0 \to Inv(G) \to A^0_G (G/S, H^*).$$
As noted in \S\ref{subsec:strat}, we have $A^0_G (G/S, H^*)=Inv(S)$.\end{proof}

\subsection{Applications}

\begin{ex}\label{ex:kleenex}
As a very simple illustration, we may compute the invariants of a
wreath product $\Gm\wr\Z/2$, that is a semi-direct product $\Gm^2
\rtimes \Z/2$ where $\Z/2$ permutes the two copies of the
multiplicative group. 

Indeed, this group may be seen as the normalizer of a maximal torus in
$\GL_2$, and thus it has a canonical representation $W$ of dimension
$2$. The space $\pr (W)$ is just two orbits, and that of $[1,1]$ is
open with stabilizer $\Gm\times \Z/2$. From corollary \ref{coro:skip}
(and corollary \ref{coro:kunneth}), we know that the restriction map
from $Inv(\Gm^2 \rtimes \Z/2)$ to $Inv(\Z/2)$ is injective. It is
also certainly surjective, since the inclusion $\Z/2 \to \Gm^2 \rtimes
\Z/2$ is a section for the projection $\Gm^2\rtimes \Z/2 \to \Z/2$.

Finally, we conclude that $Inv(\Gm^2\rtimes \Z/2)=Inv(\Z/2)$
(regardless of the choice of $p$).
\end{ex}

%\begin{ex} 
%The group $\Z/2=\{1,\tau\}$ acts on $\Z/2 \times \Z/2$ by permuting
%the factors, and the semidirect product $(\Z/2 \times \Z/2)\rtimes
%\Z/2$ is the familiar wreath product $\Z/2 \wr \Z/2$. Also, $\Z/2$ can
%be made to act on $\Gm$ via $\tau(t)= t^{-1}$. The semidirect product
%$\Gm\rtimes \Z/2$ is then $\O_2$.

%Presently we shall turn our attention to the mod $2$ invariants of
%$$G=(\Z/2 \times \Z/2 \times \Gm)\rtimes \Z/2.$$ This group will show
%up in the sequel.

%We define a representation of $G$ on $V=\aff^2$ by letting $(\epsilon_1,\epsilon_2,t)\in \Z/2 \times \Z/2 \times \Gm$ act as
%$$\left(
%\begin{array}{cc}
%\epsilon_1 t & 0
%\\ 0          & \epsilon_2 t^{-1}
%\end{array}\right),
%$$ while $\tau$ acts on $V$ by permuting the coordinates.

%Projectively, the study of the orbits is exceedingly simple, for the action of $G$ is in fact transitive on $\pr (V)$. Moreover, the stabilizer of $[1,1]$ is isomorphic to $\Z/2 \times D_8$, where $D_8=\Z/4 \rtimes \Z/2$ is the dihedral group of order $8$.

%t follows that $$Inv(G)=A^0(BG,H^*)=A^0_G(\pr (V),H^*)=Inv(\Z/2 \times D_8),$$ by lemma \ref{lem:proj} and \S\ref{subsec:strat}. We shall compute the invariants of $D_8$ when $k=\C$ in the next section.

%end{ex}

\begin{ex}\label{ex:spin7} In \cite{skip}, Garibaldi gives many examples of
  applications of corollary \ref{coro:skip}. It is possible to recover
  a good number of his results using the techniques of this paper. Let
  us illustrates this with $\Spin_7$ at $p=2$.

If $\Delta$ denotes the spin representation of $\Spin_7$, one can show
that there exists an open orbit $U$ in $\pr (\Delta)$ with stabilizer $G_2
\times \Z/2$ ($G_2$ is the split group of that type). Hence an
injection
$$0 \to Inv(\Spin_7) \to Inv(G_2\times \Z/2).$$ Garibaldi computes the
image of this, creating invariants of $\Spin_7$ by restricting
invariants of $\Spin_8$, and exploiting the fact that further
restriction from $Inv(\Spin_7)$ to $Inv(\Z/2)$ is zero, since it
factors through the invariants of a maximal torus.

Our version of this computation is to study the complement of $U$ in
$\pr (\Delta)$. It consists of a single orbit with stabilizer $K$, and there is an affine bundle map $BK\to B(\SL_3\times \Z/2)$. Therefore we have an exact sequence
$$0 \to Inv^d(\Spin_7) \to Inv^d(G_2\times \Z/2) \to Inv^{d-1}(\SL_3\times
\Z/2).$$

This seems to be typical of the stratification method when projective
representations are used: one can obtain some results very rapidly,
but they are not always as accurate as one may wish.

Instead, if we proceed exactly as in the proof of theorem
\ref{thm:exactseq}, we obtain:
$$0 \to Inv^d(\Spin_7) \to Inv^d(G_2\times \Z/2) \to Inv^{d-1}(\SL_3) \to 0.$$
\end{ex}

Since $\SL_3$ is special, $Inv(\SL_3)=H^*(k)$. From \cite{serre} we
have that $Inv(G_2)$ is a free $H^*(k)$-module on two generators of
degree $0$ and $3$. By corollary \ref{coro:kunneth}, $Inv(G_2 \times
\Z/2)$ is a free $H^*(k)$-modules with $4$ generators of degree $0$,
$1$, $3$ and $4$. It follows that $Inv(\Spin_7)$ is a free
$H^*(k)$-module on $3$ generators of degree $0$, $3$ and $4$.

\section{The Bloch \& Ogus spectral sequence}\label{sec:spec}

\subsection{The spectral sequence}

In this section, we assume that $k$ is algebraically closed. The prime $p$ being fixed as always, we shall write $H^i_\et(X)$ for the \'etale cohomology group $H^i_\et(X,\Z/p)$.

\begin{thm}[Rost, Bloch, Ogus] \label{thm:blochogus} Let $k=\bar k$, and let $X$ be equidimensional. Then there is a spectral sequence
$$E_2^{r,s}=A^r(X,H^{s-r})\Rightarrow H^{r+s}_\et(X).$$ In particular the $E_2$ page is zero under the first diagonal, and the resulting  map $$A^n(X,H^0)=CH^n X\otimes_\Z \Z/p \to H^n_\et(X)$$ is the usual cycle map.
\end{thm}

A word of explanation on the authorship of the spectral sequence. To start with, there is the well-known coniveau spectral sequence, which converges to $ H^n_\et(X)$ and for which there is a description of the $E_1$ term. In \cite{blochogus} Bloch and Ogus prove that the $E_2$ term can be identified with $H^r_{Zar}(X,\mathcal{H}^s)$, where $\mathcal{H}^s$ is the sheafification of $U \mapsto H^s_\et(U)$. They deduce that the groups are $0$ under the diagonal, and prove the statement about the cycle map. More than 20 years later, in \cite{rost}, corollary 6.5, Rost proves that the $E_2$ term can also be described using his Chow groups with coefficients as in the theorem. The sequence seems to be usually refered to as the Bloch \& Ogus spectral sequence.

Given an algebraic group $G$, we may take a model for $BG$ to play the role of $X$, and obtain:

\begin{coro} There is a map
$$H^*_\et(BG) \to Inv(G)$$ which vanishes (when $*>0$) on the image of the cycle map $$CH^*BG \otimes_\Z \Z/p \to H^{2*}_\et(BG).$$
\end{coro}

\begin{ex}
The Stiefel-Whitney classes as in \ref{subsec:ortho} can be defined as the images of the elements with the same name in the cohomology of $B\O_n$. The corollary explains why they square to zero in $Inv(\O_n)$, since they square to Chern classes in cohomology, and these come from the Chow ring (see \cite{totaro} for details).
\end{ex}

\begin{rmk} When $G$ is a finite group, viewed as a $0$-dimensional algebraic group, there is a natural map
$$H^*(G,\F_p) \to H^*_\et(BG),$$ where $H^*(G,\F_p)$ is the usual
cohomology of $G$ as a discrete group. Indeed, there is a Galois
covering $EG \to BG$ with group $G$, so we may see this map as coming
from the corresponding Hoschild-Serre spectral sequence. Since we assume in this section that $k$ is algebraically closed, it follows easily that the map is an isomorphism.

Composing this with the map from the previous corollary, we obtain the
homomorphism $$H^*(G,\F_p) \to Inv(G)$$ which was considered in
\cite{serre}. Of course the direct definition given in \lc~ is much
preferable.
\end{rmk}

\begin{coro}[to theorem \ref{thm:blochogus}] \label{coro:lowdegree} There are natural isomorphisms:
$$Inv^1(G)=H^1_\et(BG)$$ and
$$Inv^2(G)=\frac{H^2_\et(BG)}{CH^1BG\otimes_\Z \Z/p}.$$
\end{coro}

The denominator in the second isomorphism is really the image of
$CH^1BG\otimes_\Z \Z/p$ in \'etale cohomology via the cycle map. When
$k=\C$, the cycle map is injective in degree $1$ and in fact
$CH^1BG=H^2(BG,\Z)$ (topological cohomology), see \cite{totaro}.

\begin{coro}[to theorem \ref{thm:blochogus}] \label{coro:pedrorost}
Any class $x$ in the kernel of the cycle map $$CH^2 BG \otimes_\Z \Z/p
\to H^4_\et(BG)$$ determines an invariant $r_x\in Inv^3(BG)$, which is
well-defined up to the image of $H^3_\et(BG) \to Inv^3(BG)$. If $x$ is
nonzero, neither is $r_x$.
\end{coro}

We think of $r_x$ as a simplified version of the Rost invariant.

\begin{proof}
 Our assumption implies that the class $x$, viewed as an element of bidegree $(2,2)$ on the $E_2$
  page of the spectral sequence under discussion, must be hit by a
  differential. Let $r_x\in E_2^{0,3}=Inv^3(G)$ be such that
  $d_2(r_x)=x$. This element is well-defined up to the kernel of
  $d_2$. Since further differentials $d_r$ for $r>2$ are zero on
  $E_r^{0,3}$, we see that $r_x$ is defined up to elements which
  survive to the $E_\infty$ page. These are, by definition, the
  elements in the image of $H^3_\et(BG) \to Inv^3(BG)$.
\end{proof}

\subsection{Applications}

For simplicity, we shall take $k=\C$, the complex numbers, in the
applications. In this case according to remark \ref{rmk:complex}, the
\'etale cohomology of $BG$ as above coincides with the topological
cohomology of a topological classifying space (i.e. a quotient $EG/G$
of a contractible space $EG$ endowed with a free $G$-action).

We start with a proposition which should be compared with statements
31.15 and 31.20 in \cite{boi}, for which there is no proof available
as far as I am aware. 

\begin{prop}\label{prop:connected} Over the complex numbers, there are isomorphisms
$$Inv^1(G)=Hom(\pi_0(G),\Z/p)$$ and
$$Inv^2(G)=p \textrm{-torsion in}~ H^3(BG,\Z).$$ In particular, if $G$
  is connected then $Inv^1(G)=0$, and if $G$ is 
  $1$-connected then $Inv^2(G)=0$.
\end{prop}

\begin{proof} From corollary \ref{coro:lowdegree}, we have
  $Inv^1(G)=H^1(BG,\F_p)$, and of course this is
  $Hom(\pi_1(BG),\Z/p)$. However $\pi_1(BG)=\pi_0(G)$.

As noted above, the cycle map is injective in degree $1$ over the
complex numbers, and $CH^1 BG=H^2(BG,\Z)$. From corollary
\ref{coro:lowdegree} again, we see that
$$Inv^2(G)=\frac{H^2(BG,\F_p)}{H^2(BG,\Z)\otimes \Z/p},$$ and it is
elementary to show that this maps injectively onto the $p$-torsion in
$H^3(BG,\Z)$ via the Bockstein.

The statement about connected groups is trivial. When $G$ is
$1$-connected, it is automatically $2$-connected, since any real Lie
group has $\pi_2(G)=0$. Thus $BG$ is $3$-connected and we draw
$H^3(BG,\Z)=0$ from Hurewicz's theorem.
\end{proof}

\begin{ex} Consider the case of the exceptional group $G_2$, for
  $k=\C$ and $p=2$. It is reductive and $1$-connected, so by the
  proposition we have $Inv^i(G_2)=0$ for $i=1,2$. Moreover, the Chow
  ring of $G_2$ has been computed over the complex numbers, see
  \cite{pedrog2}. It turns out that the map $CH^2BG_2 \otimes_\Z \Z/2
  \to H^4(BG_2,\F_2)$ has exactly one nonzero element. From corollary
  \ref{coro:pedrorost}, we know that there is a nonzero invariant
  $e_3\in Inv^3(G_2)$. It is uniquely defined as $H^3(BG_2,\F_2)=0$.

It is proved in \cite{serre} that for any field $k$ of characteristic
$\ne 2$, $Inv(G_2)$ is in fact a free $H^*(k)$-module on the
generators $1$ and $e_3$.
\end{ex}

\begin{ex}\label{ex:d8}
We take $k=\C$, $p=2$, and $G=D_8=\Z/4\rtimes \Z/2$, the dihedral group. We shall completely compute $Inv(G)$ by showing first that things reduce to corollary \ref{coro:lowdegree}.

Let $t$ be a generator of $\Z/4$ in $G$, and let $\tau$ be the second generator, so that $\tau(t)=\tau t \tau^{-1}=t^{-1}$. We let $G$ act on $V=\aff^2$ via
$$t\mapsto \left( \begin{array}{cc} t & 0 \\ 0 & t^{-1} \end{array}\right)$$
and $$\tau \mapsto \left( \begin{array}{cc} 0 & 1 \\ 1 & 0 \end{array}\right).$$

The corresponding map $G\to\GL_2$ is an embedding, so $G$ acts freely outside of a finite number of lines. Let $U$ denote the open complement. We have an injection 
$$0 \to A^0_G(V,H^*)=Inv(G) \to A^0_G(U,H^*)=A^0(U/G,H^0),$$ by example \ref{ex:freeactions}. Now, $U/G$ is a variety of dimension $2$ over an algebraically closed field, so $H^*(k(U))=0$ for $*>2$. It follows that $Inv(G)$ is concentrated in degrees $\le 2$.
 
Since we are working over $\C$, we have $H^*_\et(BG)=H^*(BG,\F_2)=H^*(G,\F_2)$. The cohomology of $G$ is well-known, see for example \cite{adem}:
$$H^*(G,\F_2)=\F_2[x_1,y_1,w_2]/(x_1y_1=0).$$ More precisely, $G$ can be presented as an "extraspecial group", ie as a central extension
$$0 \to \Z/2 \to G \to E \to 0$$ where $E\approx \Z/2 \times \Z/2$. If we choose $x_1$ and $y_1$ so that $H^*(E,\F_2)=\F_2[x_1,y_1]$, then we can pullback these classes to the cohomology of $G$ where they will give the classes with the same name in the description above. (Besides, the cohomology class of the extention is $x_1y_1$.)

We have immediately, by corollary \ref{coro:lowdegree}, that $Inv^1(G)=\F_2\cdot x_1 \oplus \F_2\cdot y_1$. In degree $2$, we note that $H^2(G,\F_2)$ is generated additively by the classes $x_1^2$, $y_1^2$, and $w_2$. The Chow ring of $BE$ is $CH^*BE \otimes_\Z \Z/2=\F_2[x_1^2,y_1^2]$, from which we know that the classes $x_1^2$ and $y_1^2$ in the cohomology of $BG$ certainly come from the Chow ring of $BG$.

On the other hand, $Sq^1 w_2=w_2(x_1 + y_1)\ne 0$ (\lc), while the Steenrod
operation $Sq^1$ is zero on classes coming from the integral
cohomology, and {\em a fortiori} it is zero on the classes coming from
the Chow ring.

As a result, we have finally $Inv^2(G)=\F_2\cdot w_2$.

Over a general field, we could reach a similar conclusion by studying the geometric situation a bit more carefully. The variety $U/G$, for example, can be shown to be the open subset in $ \aff^2$ obtained by removing the axis $Y=0$ and the two parabolae $X\pm 2Y^2=0$. Alternatively, you might want to use that $D_8$ is a $2$-Sylow in $S_4$, and exploit the double coset formula as in \cite{serre}, chap. V.

\end{ex}

\section{Other cycle modules}\label{sec:milnor}

We shall conclude with a few simple remarks on other possible
cycle modules.

Given any cycle module $M$ as in \S\ref{sec:rost}, we may define the
invariants $Inv(G,M)$ of $G$ with values in $M$ to be the natural
transformations of functors
$$H^1(-,G) \to M(-).$$ It is straightforward to establish the
inclusion $$Inv(G,M)\subset A^0(BG,M)$$ for any $BG$ which is the base
of a versal $G$-principal bundle. Indeed, the arguments in
\cite{serre} hold {\em verbatim} (and in fact in \lc~ the reader will
find a similar inclusion even for invariants with values in the Witt
ring, even though the Witt ring satisfies weaker properties than cycle
modules do.)

We may also define the equivariant Chow groups $A^*_G(X,M)$ for
varieties $X$ acted on by $G$, exactly as we have done for
$M=H^*$. When this is done, we may take $BG$ to be $\spec (k)\bor G$
(as in the rest of this paper), and we have an equality
$$Inv(G,M)=A^0(BG,M)=A^0_G(\spec (k),M)$$ by Totaro's arguments as in \cite{serre}.

The techniques we have used for Galois cohomology may be used for any
cycle module. Let us illustrate this for $M=K_*$, the algebraic
$K$-theory of fields (Milnor or Quillen, it will not affect the
sequel). In this case we have $A^*(X,K_0)=CH^*X$, and the parallel
with our previous computations becomes even more obvious.

Arguing as in example \ref{ex:gm}, we obtain that $A^0(\Gm,K_*)$ is a
free $K_*$-module on two generators, one in dimension $0$, the other
in dimension $1$. If $p$ is a prime number, we may use this to compute
the invariants of $\Mu_p$.

The group $\Mu_p$ has a $1$-dimensional representation $V$, and $\Gm$
sits in $V$ as a $\Mu_p$-invariant open subset whose complement is a
point $\spec (k)$. We obtain the exact sequence (where $G=\Mu_p$):
$$0 \to A^0_G(V,K_*) \to A^0_G(\Gm,K_*) \to A^0_G(\spec (k),K_{*-1})
\to A^1_G(V,K_{*-1}).$$

The action on $\Gm$ is free with quotient $\Gm$, so as in example
\ref{ex:freeactions} we draw $A^0_G(\Gm,K_*)=A^0(\Gm,K_*)$. Thus we
may rewrite the exact sequence:
$$\begin{CD}
0 @>>> Inv(G,K_*) @>>> A^0(\Gm,K_*) @>>> Inv(G,K_{*-1}),
\end{CD}$$
and for $*=1$ it is worth writing the extra term:
$$0 \to Inv(G,K_1) \to A^0(\Gm,K_1) \to Inv(G,K_{0}) \stackrel{s_*}{\longrightarrow} CH^1BG \to 0.$$
Let us explain the notation $s_*$ and what this map looks like. We
call $s$ the zero section $s:BG \to V\bor G$ of the vector bundle
$\pi: V\bor G \to BG$, and $s_*$ is the induced pushforward map. We
know that $\pi^*$ is an isomorphism, and that there is a projection
formula $s_*(\pi^*(x)y)=xs_*(y)$. If we use $\pi^*$ as an
identification, it follows that $s_*(x)=c_1(V)x$ where $c_1(V)=s_*(1)$
is the first {\em Chern class} of $V$.

Finally the map $s_*: Inv(G,K_0) \to CH^1 BG$ is simply the surjective
map $\Z \to CH^1BG$ sending $1$ to $c_1(V)$. Now, $CH^1BG$ is
$p$-torsion by a transfer argument (or see \cite{totaro}, example
13.1, which shows that $CH^1 BG$ may well be $0$ depending on $k$). In
any case, $s_*$ has a kernel isomorphic to $\Z$.

We conclude that the $1$-dimensional generator for $A^0(\Gm,K_*)$ must
map to a generator for this kernel. As a result, $Inv(G,K_*)$ is
reduced to $K_*$, i.e. {\em the group $\Mu_p$ has no nonconstant invariants in
  algebraic $K$-theory at all}.

For $p=2$ for example, assuming that $char(k)\ne 2$, we may use the
surjection $H^1(K, (\Z/2)^n) \to H^1(K,\O_n)$ for fields containing
$k$ to deduce that $Inv(\O_n,K_*)$ injects into
$Inv((\Z/2)^n,K_*)$. Thus, after an obvious K\"unneth argument, we see
that $\O_n$ has no nonconstant invariants in algebraic $K$-theory,
either. In particular, there is no natural way of lifting the
Siefel-Whitney classes to integral Milnor $K$-theory.

\bibliography{myrefs}
\bibliographystyle{siam}
\end{document}